\documentclass[preprints,mathematics,
article,submit,sort,
moreauthors,pdftex]{Definitions/mdpi} 
\usepackage{cjhebrew}
\usepackage{comment}
\firstpage{1} 
\pubvolume{xx}
\issuenum{1}
\articlenumber{5}
\pubyear{2019}
\copyrightyear{2019}
\externaleditor{Academic Editor: Vladimir V. Mazalov}
\history{Received: date; Accepted: date; Published: date}

\usepackage{filecontents}
\usepackage{subcaption}

\usepackage{diagbox}
\usepackage{bm}
\graphicspath{{Figures/},{Pictures/}}
\RequirePackage[hyperpageref]{backref} 
    \renewcommand*{\backref}[1]{}  
    \renewcommand*{\backrefalt}[4]{
       \ifcase #1 
          No cited.
       \or
          Cited on p. #2.
       \else
          Cited on pp. #2.
       \fi}

\newcommand{\nthmname}{Theorem}\newcommand{\exaname}{Example}
\newcommand{\ndefname}{Definition}\newcommand{\nremname}{Problem}\newcommand{\ncorname}{Corollary} 
\usepackage{amsfonts}

\DeclareMathOperator{\bbP}{\bold P}

\preto{\abstractkeywords}{\nolinenumbers}

\Title{Drivers' skills and behavior vs. traffic at intersections}

\Author{Krzysztof J. Szajowski $^{1,\dagger,\ddagger}$\orcidA{} and Kinga Włodarczyk $^{2,\ddagger}$\orcidB{}}

\AuthorNames{Krzysztof Szajowski, Kinga Włodarczyk}

\address{%
$^{1}$ \quad Wrocław University of Science and Technology; Krzysztof.Szajowski@pwr.edu.pl \\
$^{2}$ \quad Wrocław University of Science and Technology;; kingawlodarczyk96@gmail.com}

\corres{Correspondence: Krzysztof.Szajowski@pwr.edu.pl; Tel.:  +48-603-765-541 (K.Sz.)}

\firstnote{Current address: Wrocław University of Science and Technology, Faculty of Pure and Applied Mathematics, Wybrzeże Wyspiańskiego 27, PL-50370 Wrocław} 
\secondnote{These authors contributed equally to this work.}

\abstract{The work is devoted to ways of modeling street traffic on a street layout without traffic lights of an established topology. The behavior of traffic participants takes into account the individual inclinations of drivers to creatively interpret traffic rules. Participant interactions describe game theory models that provide information for simulation algorithms based on cellular automata. Driver diversification comes down to two types often considered in such research: {\bf DE}(fective)-agent and {\bf CO}(operative)-agent. Various ways of using the description of traffic participants to examine the impact of behavior on street traffic dynamics were shown. Directions for the further detailed analysis were indicated, which requires basic research in the field of game theory models.
}



\renewcommand\labelenumi{(\roman{enumi})}
\renewcommand\theenumi\labelenumi

\begin{document}
\tableofcontents
\listoffigures
\listoftables

\section{Introduction.} The considerations of this work relate to traffic modeling at intersections or places where vehicle interaction can change the behavior of a vehicle stream. The area of the problem formulated in this way is relatively wide due to the possibility of examining both complex situations and individual events on the road. Although objects occurring in traffic are of considerable size, similarity to mesoscopic systems is observed, i.e. those that fit between the micro world (in systems of single atoms or molecules described by quantum mechanics) and the macro world (in objects consisting of a very large number of particles, subject to the laws of classical mechanics). A comprehensive review of the literature on various issues of modeling circular motion can be found in the paper by Albi et al.~\cite{Alb2019:Vehicular}. We are interested in linking statistical behavior of drivers with the dynamics of vehicle streams. To this end, we must select the appropriate model describing the behavior of the driver-vehicle system and link it to the description of the stream of vehicles on the roads in the analyzed area. Although probabilistic modeling of driver behavior and vehicle traffic is a natural research method (v. \cite{UseCenGom2017:Work}, \cite{ChaChow2014:Stat}), the combination of these two elements constantly leaves many questions and is the subject of research. Calibration methods are difficult, we only have access to many factors describing the vehicle-driver system through symptoms. That is why we are looking for a link between observable driver characteristics and vehicle stream parameters. In these considerations, mathematical modeling of behavior using game theory methods should become a facilitation.

A review of the literature related to traffic modeling at intersections is very extensive, and vehicle traffic on the road is considered in many aspects. Research on theory and modeling of traffic at intersections began in the 1930s.The paper \cite{Gre1935:Traffic} by Greenshields was a pioneering work in this field. He used photographic measurement methods to calculate traffic volumes, predict and explain possible observation trends in reality. The interest in this field has increased significantly since the nineties, mainly due to the high development of road traffic. To study the dynamics of motion, many models have been proposed, e.g. Zhang's   hydrodynamic models proposed in \cite{Zha1998:nonequilibrium} and a list of scientific journals in 1994--2002, or models based on gas kinetics developing in a similar period, e.g. proposed in 2001 by Helbing et al.~\cite{gass2001:Master}. A different approach was presented by models related to car models, and the breakthrough was the use of cellular automata in 1992, which proved to be an easy and efficient method of modeling movement. Several models have been developed to depict various aspects of road traffic. At the same time, an area was developed dealing with issues related to the interaction of drivers at intersections. One approach in modeling driver behavior at an intersection refers to concepts related to game theory (cf. \cite{LiuXinAdaBan2007:Game}, \cite{KitKei2002:Game}). There were also many publications focused on the proper selection of traffic lights so as to minimize the number of cars waiting before the intersection, and research was conducted on the psychology of human behavior itself, which was also reflected in modeling of traffic.

We will focus on the impact of the behavior of individual drivers in moments of interaction with another driver on the ownership of the stream of vehicles (cf. \cite{YaoJiaZhoLi2018:BestResponse}). Cellular automata are mathematical objects for modeling various phenomena (cf.  \cite{Ila2001:Automata}). 
{The creator of cellular automata is Janos von Neumann~\cite{Neu1966:Self}, a Hungarian scientist working at Princeton. In addition, Lviv mathematician Stanisław Ulam~\cite{Embryo2010}, who was responsible for the discretization of time and space of automata, and considered the creator of the term cellular automata as "imaginary physics" had a significant impact on the development of this area.}They are used in many areas of science, in modeling physical phenomena and interactions between objects. They are also applied in modeling of movement.  The basic knowledge of the traffic simulation method used in the paper can be found in the paper by Małecki and Szmajdziński~\cite{MalSzm2013:AnalizaRuchu} (cf. also \cite[sec. 2]{BerTag2017:CAmath},  \cite{Zyg2019:AuKo}).  According to \cite{Ila2001:Automata}, cellular automata can reliably reflect many complex phenomena using simple rules and local interactions. They are a network of identical cells, each of these can assume one specific state, with the number of states being arbitrarily large and finite. The processes of changing the state of cells are carried out parallelly and in accordance with applicable rules, usually depending on the current state of the cell or the state of neighboring cells. Basically, three ways of modeling traffic at intersections described by driver interactions were selected in the research presented in this paper. Assumptions, resulting rules and evaluation of consequences for participants were introduced into the functioning of the sample intersection network. Each of the models is simulated using cellular automata. 

Because in reality, drivers do not always follow the rules of traffic, their behavior will be linked to with traffic parameters in the constructed and analyzed models. It has been assumed that drivers generally comply with the provisions of the Highway Code, but in some situations they tend to depart from them and break the rules, causing disruptions resulting in a slowdown of traffic (cf. \cite{TanKukHag2014:Social}). Details of driver behavior modeling are included throughout the work, and a more comprehensive introduction is provided in Section~\ref{KWDriDecMod}. Although there are more possible models in this area, which we will mention later, we will focus on three of them. The precise description is given below in Section \ref{KWKSzManOrg}. The considered drivers' interactions are modeled using the game theory apparatus and methods of the mass service theory.

In traffic modeling or the use of transport, game theory methods appear naturally. In the models analyzed in this work, the strategic behavior of drivers is of an auxiliary nature, which will be signaled in the right places when introducing and analyzing models. For a more complete picture, in the next section, we will signal some other traffic problems analyzed by creating mathematical models.

\subsection{\label{KWDriDecMod}Driver decision models.} 
		\textls[-15]{Fisk \cite{Fis1984:Game} in~1984 described correspondences between two game theory models (Nash noncooperative and Stackelberg games) and some problems in transportation systems modeling. An example of each is described in detail, namely the problem of carriers competing for intercity passenger travel and the signal optimization problem. The discussion serves to underline differences between two categories of transportation problems and introduces the game theory literature as a potential source of solution algorithms. Also, it is shown that inner-outer iterative techniques for Stackelberg type problems cannot be expected to converge to the solution, and an approximate formulation of these problems is introduced which appears to be more readily solvable. However, this discussion is far away from determining driver modeling. Here, the two equilibrium concepts, Nash and Stackelberg respectively, can be used to discover which action set or strategy is optimal for every participant in the game. The participants are drivers. Optimality in this context is evaluated on the basis of payoffs resulting from the decisions by (and interaction among) the participants (v. the monographs by Ferguson \cite{Fer1967:MS}, Owen~\cite{Own2013:Game}, Platkowski~ \cite{Pla2012:Introduction}, Mazalov~\cite{Maz2014:book}. Payments in traffic modeling games come down to passing times -- their shortening or lengthening.} 
	
In issues that have the common trait that decision-makers know that their result cannot be achieved at the expense of the other community, one cannot rely solely on pure antagonist game models. Instead of talking about modeling the game, it's better to think about modeling the behavior of project participants. The existing objective dependencies mean that decision-makers are motivated to take into account these dependencies and generally do not act independently, although they are not able to agree their actions and form formal coalitions. One can only assume that they are motivated to coordinate their proceedings. This, in turn, forces us in modeling to adopt appropriate sets of strategies or otherwise model information available to players. Achieving a common optimal result in an orthodox model of game theory does not introduce a general reason or justification for choosing the right strategies. It is known that in the simplest cases, participants in a joint project generally easily coordinate their decisions without difficulty. The recognition of this in the mathematical model is not known today, because the actual mechanisms of such coordination, the way to achieve it is poorly understood. There are theories explaining strategic coordination, but their implementation in the mathematical model has limited application. The reason for this is the need to change the specifications of the game and make incredible assumptions. By adopting Stackelberg's extreme rationality, according to which players only choose strategies that maximize their own profits, in conditions where co-decision makers can always foresee opponents' strategies and respond to them as best as possible, I avoid these problems. This makes it possible to clarify strategic coordination in the common interest of all project participants. Previous experimental encouraged this approach. They showed that Stackelberg's approach in asymmetrical games is rational.
	
Only vehicles are regarded as the game participants. Kita et al.  \cite{KitKei2002:Game} has adopted a game theoretic analysis to consider a merging-give way interaction between a through car and a merging car, which is modeled as a two-person non-zero sum non-cooperative game. Kita’s approach can be regarded as a game theoretic interpretation of Hidas' driver courtesy considered in \cite{Hid2002:Merging} from the viewpoint that the vehicles share the payoffs or heuristics on the lane changes, which is a reasonable traffic model but fails to assign uncertainties resulting from the action of the other human drivers. Moreover, one cannot guarantee that the counterpart would act as determined in the game since the counterpart may be able to consider other factors that the subject driver cannot take into account. Accordingly, it is necessary to design an individual driver model that does not share their payoffs in the decision making processes to reflect such an uncertainty. This approach, as we shall see later, facilitates a more realistic model of driver behavior in traffic situations. It can be found e.g. in \cite{LiuXinAdaBan2007:Game}.

The behavior of drivers crossing the intersection or joining traffic from another road is a potential source of conflict with another road user. An additional element intensifying the conflict are various assessments of the situation resulting from different levels of skills and the ability to use them. The controllers (drivers) can be divided into two types in a simplified way  (cf.  Paissan and Abramson \cite{PasGui2013:Imitation}, Fan et al. \cite{FanJiaTiaYun2014:TraficGT}, and Yao et al. \cite{YaoJiaZhoLi2018:BestResponse}):
\begin{itemize}\label{CODEdrivers}
\item \textbf{Regulatory drivers}. They will be abbreviated as \textbf{CO} (\textit{Cooperator}).
\item \textbf{Non-compliant drivers}. They will be abbreviated as \textbf{DE} (\textit {Defector}).
\end{itemize} 
Although we see a natural possibility of distinguishing between a lack of knowledge of traffic rules and their non-application (conscious or unconscious), we leave such detailed analyzes for further research. In the subsequent analysis assuming that each of them will react according to the category.

Traffic models and driver behaviors are generally generic and require calibration to suit their place and time of use. Performing such a procedure requires obtaining relevant data and the use of adequate statistical methods. The specificity of the problem leads to the formulation of basic research problems in both modeling and statistics. This is signaled by numerous publications on road engineering, road safety and related problems, such as driver behavior (cf.DIng and Huang~ \cite{DinHua2009:TrafficFlow}, Bifulco et al.~ \cite{BifGalParSpeGai2014:TrafficData}). Some questions may be solved by choosing and adapting models known as decision theory. One of the elements discussed in this article are the skills and behavior of drivers. We suggest using Bayesian and minimax estimation methods to assess the parameters associated with modeling the distribution of drivers' characteristics (cf. \cite[p. 17]{BicDok2015:MS},  \cite{Ber1980:SDT}).

\subsection{\label{KWinterdritra}Intersections, drivers and traffic.} Intersections are an inseparable element of road traffic (cf.\cite[Section 1A.13, def. 94]{MUTCD2009}). In this consideration it is assumed that they are equal without junctions regulatory
({for further information concerning classification of crossroads can be found e.g. at the OSK Duet driving school website, \emph{Virtual driving school}, \url{http://oskduet.pl}}). At such intersections, priority is given to road signs defining one of the roads as the main road and the other as the subordinate one. In the absence of signs, the so-called right-hand rule that gives priority to all vehicles on the road on the right. Not all road users obey the rules cited above (cf. the dichotomous classification of the drivers above at page \pageref{CODEdrivers}). It often happens that drivers enforce the right of way at intersections, thus forcing other traffic participants to slow down, or sometimes causing collisions or traffic accidents. The effects of such behavior will be further explored in the work presented. When developing the research that is the subject of this work, it is worth remembering that the intersections are different and you need to consider the topology of the intersection in mathematical models.

In order to create a model of vehicle movement, we will distinguish a description of the behavior of individual participants (vehicle -- driver) and a description of the dynamics of the location of all vehicles in the analyzed region. We describe the changes in the position of the vehicle in the intersection using the cellular automata method and Nagel and Schreckenberg's model (\textbf{NaSch}, v. \cite{NagSch1992:freeway}) described in the section \ref{modelnagela}. This is a proven method that allows testing the impact of changes in driver behavior on vehicle flow parameters. However, in the description of the dynamics of traffic at the intersection, three elements can be distinguished, the specification of which is important for aspects of interest to us. Those are:
\begin{enumerate}
\item \textbf{Identifying road elements} boils down to the rules of right-hand traffic. This means that priority is given to the one on the main road when meeting at the intersection of two vehicles, i.e. the one who sees the second vehicle on its left (priority of the road on the right-hand side). This objective determination is transformed by decision makers. It is known that the main reason for perturbation in the stream of vehicles are driving behaviors that do not comply with traffic rules (v. \cite{TanKukHag2014:Social}, \cite{Cor-BerGerSte2016:Traffic}). Earlier studies by Mesterson-Gibbons~\cite{MesGub1990:Dilemma} have found various quantification of driver behavior, however, two categories of drivers have been adopted for the purposes of this study. By convention, they will be those who follow the rules of the road (\textbf{CO} drivers) and those who do not comply with them (\textbf{DE} drivers). With the approach used to analyze the phenomenon, the proportions of the types of drivers present in the population in the analyzed area are significant.
\item \textbf{Determining the behavior of drivers} is based on the fact that, at each meeting, before the settlement of traffic in the next step  (nearest second), the type of driver is identified, and, on this basis, his decision is determined, which translates into vehicle behavior. There are several ways to identify the types of individual drivers in the considerations. In the models selected for detailed analysis, the method of determining the driver's behavior is different and depends on the assumptions made earlier in relation to the rules functioning in everyday life.
\item \textbf{Priority assignment} at the intersection results from setting their priorities resulting from the types described above and assigned to drivers. Costs (in units of increasing or decreasing speed) related to strategies adopted by drivers were determined. Depending on the adopted model, the payout values are different.
\end{enumerate}

\begin{figure}[tbh!]
\centering
\includegraphics[width = 8cm]{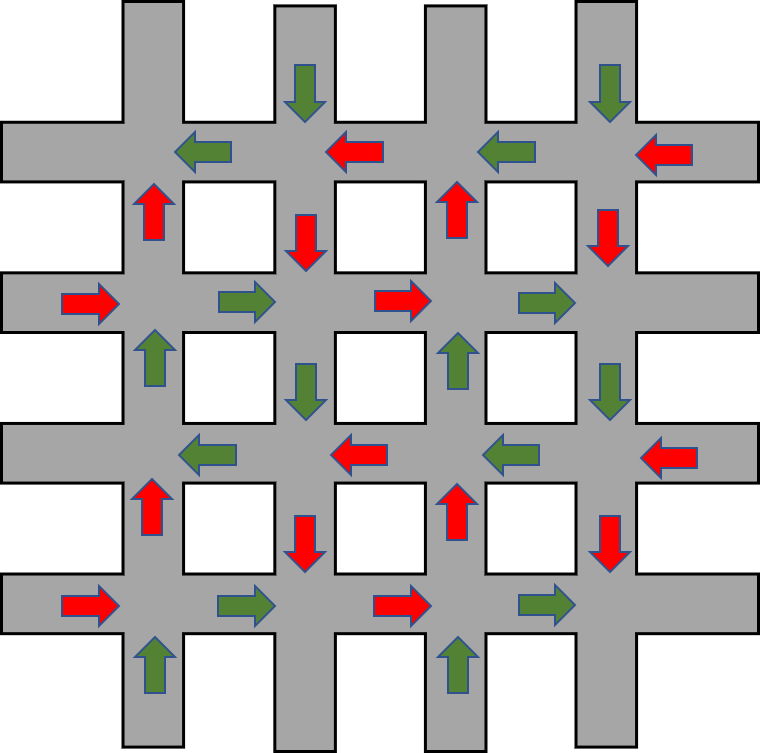}
\caption{\label{KWnetwork} The street system under consideration. } 
\end{figure}

For our research, we accept street topography previously used in \cite{PasGui2013:Imitation} or \cite{MesGub1990:Dilemma}. The traffic system consists of a network of equivalent streets, 4 of which run north-south and 4 east-west, forming a regular grid. Each road is single-lane and one-way, but the directions of vehicle traffic are different. Cars on two of the horizontal streets move from right to left, and on the other two in opposite directions. The same situation occurs in the case of streets arranged vertically, the direction of two is facing downwards, and the others are facing upwards. There are no right-of-way streets, so the right hand rule applies in the presented system. 16 road junctions and four directions of travel are possible: from top to bottom, from bottom to top, from right to left, from left to right. Therefore, four types of possible meetings of drivers on the road were received and each other driver gets the right of way. In Figure \ref{KWnetwork}, vehicles and their direction are shown by arrows. The arrows indicate vehicles and their return indicates the direction of travel. Green symbolizes the right of way, and red means that the car, according to the right hand rule, gives way. This topography is a simplification. Thanks to symmetries and uniformity, traffic analysis at such an intersection is easier and allows for proper interpretation of the results.

\subsection{\label{modelnagela}Nagel-Schreckenberg model.} 
An efficient method of traffic simulation was proposed by Nagel and Schreckenberg in 1992 by German physicists, published in  \cite{NagSch1992:freeway}. It presents the movement of cars on a straight one-lane and one-way road. The road was divided into 7.5-meter sections corresponding to the average length of an average car along with the distance in front and behind the car. Each of these sections is represented by a single cell of the automaton. The cell can be empty or occupied by one vehicle. Each vehicle $ i $ has a specific speed $ v_i $, which informs about the number of cells it will travel in one time step, with the speed not exceeding the set maximum speed $ v_{max} $. The transition function (v. \cite[sec. 2]{BerTag2017:CAmath},  \cite{Zyg2019:AuKo}) responsible for the movement of vehicles consists of 3 stages, occurring simultaneously for all objects:
\begin{enumerate}
\item \textbf{Acceleration/ Braking.} The car increases its speed by one, if it is not higher than the maximum speed and the number of free cells in front of it. When the distance to the car ahead is less than the current speed, the vehicle slows down to a value equal to the empty space in front of it. In the mathematical notation it looks like this
\begin{equation}\label{przypieszenie}
v_i (t + 1) = \min (v_i (t) + 1, v_{max}, d_i).
\end{equation}
\item \textbf{Random event.} A car with a certain probability decreases its speed by 1, provided it is not less than zero. The equation for the described situation  is as follows
\begin{equation}\label{randomness}
v_i (t + 1) = \left\{
					\begin{array} {ll}
						\max (v_i (t + 1) -1.0), & \textrm{with probab. $ p $} \\
						v_i (t + 1), & \textrm {with probab. $ 1-p $.} \\
					\end{array} \right.
\end{equation}
\item \textbf{Update position.} The car moves as many cells as its current speed. According to the formula
\begin{equation}\label{update}
		x_i (t + 1) = x_i (t) + v_i (t + 1).
\end{equation}
\end{enumerate}
The \textbf{NaSch} model reliably reflects the movement of vehicles on the road and the mutual interactions of drivers. One example that is noticed when analyzing the results of simulations is the occurrence of start-stop waves. Showing how sudden braking of one driver affects other road users.

\begin{figure}[th!]
{\noindent
\small
\begin{subfigure}[b]{0.48\textwidth}
\centering {\includegraphics[width=\textwidth]{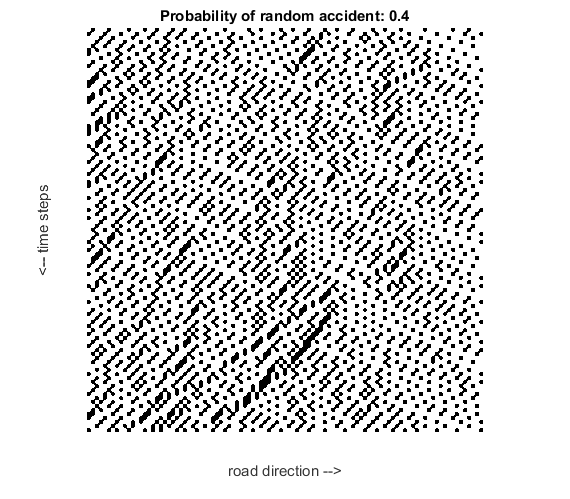}}
\caption{\label{KWnasch1a}New object appears with probability $0.4$.}
\end{subfigure}
\hfill\begin{subfigure}[b]{0.48\textwidth}
\centering {\includegraphics[width=\textwidth]{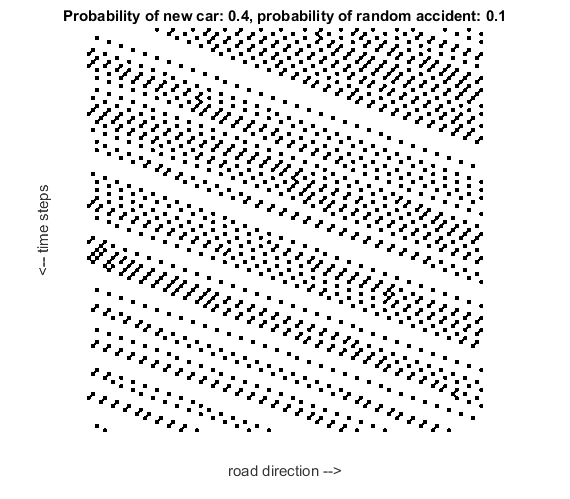}}
\caption{\label{KWnasch2a}New object appears with probability $0.1$.}
\end{subfigure}
}
\vspace{-1mm}
\caption{\label{KWnasch1i2}Simulation of motion according to the assumptions of the \textbf{NaSch} model.}
\end{figure}
Figure \ref{KWnasch1i2} presents a cellular automaton built in accordance with the assumptions of the Nagel-Schreckenberg model.  Black cells symbolize the presence of a car, and white symbolizes its absence. The simulation was performed for a periodic boundary condition. The vehicle leaving the space of the machine goes to the queue of cars waiting to enter the network, where the appearance of a new object occurs with a certain probability, which is $0.4$ at the experiment presented at Figure~\ref{KWnasch1a}. Cars move to the right, each subsequent row in the presented figure illustrates the next iteration of the program -- traffic situations in the next time step. In places, there are temporary densities of cars, caused by a sudden stop of one of the drivers, as a result of which others are also forced to slow down. The occurrence of higher densities is called backward plug. This phenomenon is observed as the effect of traffic lights. The presented situation corresponds to a small density of congestion, because one car releases one car behind it, then returns to traffic. It can be seen that stopping one vehicle causes a chain of stops for subsequent cars. For the comparison, Figure~\ref{KWnasch2a} presents the situation for lower traffic density and lower probability of random events, where a smoother process of vehicle movement is noticeable.

\subsection{\label{KWKSzManOrg}Manuscript organization.} The purpose of the research presented here is to analyze the impact of behavior, in particular interactions, of predefined types of drivers on traffic performance. Three ways of drivers' type influence are modeled and their behavior consequence on traffic at intersections are investigated. In each of the analyzed problems, we examine different aspects of traffic at the intersection. The first presented model researches the case of the constant probability of individual types of drivers. It shows how the presence of non-compliant drivers negatively affects the quality of road traffic. When costs received by drivers in conflict situations are quite high, no collisions occur frequently and there are not so many non-compliant drivers on the roads, and most road users are aware of the consequences of reckless driving. This is disused in section~\ref{DetAnalysisModel1}. 
In the next section, we focus on modeling the psychological aspects of road participation. This is an extension of the research of the section~\ref{DetAnalysisModel1} in the sense that we pay more attention to the behavior of participants due to their tendency to violate traffic law and to cooperate with other road users. Consequently, the research of the section~\ref{model2} shows that, given a certain group of drivers resistant to imitation strategy and always deciding to comply with the law, we are able to influence the final distribution of types of behavior.
The model analyzed in the section~\ref{model3} is a special case of that in the section~\ref{model2}. However, we examine here the reasons why drivers violate traffic rules as a link to increased traffic congestion. As a consequence, we believe that traffic disorder increase improper behavior. It has been detected that above a certain degree of congestion the traffic situation reaches certain limits of good performance. 
Each of the models is simulated using cellular automata. A summary of the considerations in the sections~\ref{DetAnalysisModel1}--\ref{model3} is contained in the section \ref{ModComp}. The proposed modeling of street traffic allows the study of real traffic and, as a consequence, the determination of parameters not known a priori, such as the participation in the traffic of non-compliant drivers, delay time or percentage speed delay, which significantly increases the tendency to behavior causing further problems in road traffic. This aspect is the subject of the section~\ref{conclude}.

\section{\label{DetAnalysisModel1}Simulation analysis and discussion of generalizations for Model I. }
\unskip
\subsection{Model description.} The purpose of the first model is to check the impact of drivers who do not comply with traffic rules on its overall functioning. The driver type is generated with a certain probability. The chance to draw a driver who complies with the rules is $ p $, while for a driver who does not comply with traffic rights, $ 1-p $. Paissan and Abramson in  \cite{PasGui2013:Imitation} introduce a periodic boundary condition. After leaving the network, the cars are placed in a queue, from which they go to a randomly chosen road, regardless of the street they left. Traffic is updated in accordance with the assumptions of the Nagel-Schreckenberg model. The proposed model assumes a maximum speed of 1. This assumption does not correspond to the real performance of vehicles on the road, but makes it possible to explain the suggested payouts in games. The games imitate the meetings of drivers at intersections, and the payout matrix informs about the costs incurred during the \cite{Maz2014:book} maneuver. During the meeting of drivers at the intersection, \textbf{4 scenarios} are possible:
\begin{enumerate}
\item Both drivers are cooperators (CO).
\item The driver driving on the right, i.e. the one with priority is a non-compliant driver.
\item The driver on the left, i.e. the one who should give way, is a driver who does not comply with the rules.
\item Both drivers do not follow traffic rules.
\end{enumerate}
In the first situation, both drivers comply with the rules, so one of the drivers will give way. The second scenario will not end in a conflict either, because the driver who is about to give way complies with the rules. The other two options do not have a definite solution. In the third event, the driver with the priority does not give up the road. It is assumed that this will result in the loss of $ a $ time for each participant. The most conflicting is the last case when two leading non-compliant drivers meet at the intersection, and therefore there is a risk of collision, which is a waste of time greater than before. The time the driver needs to cross the intersection is used as payoff values in games represented as non-zero matrix games (v. 
\cite{HarSel1988:GTgames,Pla2012:Introduction,Maz2014:book}). An example is Table \ref{tab1}, showing the time it takes to cross the intersection in units of simulation steps.
\begin{table}[H]
		\begin{center} 
			\caption{\label{tab1}Costs of an interaction between different type of drivers.
			}
			\vspace{0.2cm}
			\begin{tabular}{|c|c|c|} 
				\hline
				{\backslashbox{Left\strut}{\strut Right}}  & CO & DE\\
				\hline
				&&\\
				CO & 2,\; 1 & 2,\; 1\\[2ex]
				\hline
				&&\\
				DE & $d_{DC}^{ld},\; d_{DC}^{rd}$ & $d_{DD}^{ld},\; d_{DD}^{rd}$\\[2ex]
				\hline
			\end{tabular}
		\end{center}
	\end{table}
The left driver is approaching the intersection on the left, so he has no the right of the way and he has no priority. He suffers a loss of 2, these are two simulation steps he needs to take to give way and then cross the intersection. Parameters(payoffs) could be $d_{DC}^{\cdot d}=a$ and $d_{DD}^{\cdot d}=b$ to show the costs incurred in the event of a collision on the road. In~  \cite{PasGui2013:Imitation} it is emphasized that the adoption of the same costs for the left driver \textbf{DE} and the right \textbf{CO} is a simplification. In fact, it's more complex and usually \textbf{CO} delay more than \textbf{DE}. As a possible extension, they suggest $d_{DC}^{r d}=c>a+1$. The further discussion of the issue will be given in section \ref{KWKSzExtModInt}, where the type of driver is interpreted in terms of their strategies. 

The above situation corresponds to an event in which both drivers are punished for the conflict, but to avoid a collision, one waits longer. It would be necessary to agree who will be considered submissive. In the simulation below the originally proposed values were retained. There will be a situation when the drivers at the intersection "overlap". However, due to the fact of earlier expectations, it is assumed that in real life a collision at the intersection would not occur. One of the drivers would allow an opponent to cross over.

In the analyzed case, described above, it was assumed that the occurring drivers belong to two categories and each driver belongs to one of them. The ratio of driver types is constant and known. This allows simulation testing of the consequences of such an assumption. The described model can also be used to analyze real traffic (real driver behavior) to determine this ratio, as in  section~\ref{KWKSzExtModInt}.

\subsection{\label{symulacjam1}Simulations.} A network of streets with a length of $50$ cells will be considered. In the article \cite{PasGui2013:Imitation} it was proposed to introduce penalties for a collision of $100$, but it was reduced and the value $50$ was assumed at work. The cost incurred in the conflict: the left driver who did not comply with the rules, the right driver who followed the rules in accordance with the proposal used the values of $3$. In the first significant step, a driver queue was created and the correct types were assigned to them with certain probability. The queue reaches the maximum number of vehicles on the road. Vehicles presented in it are introduced into the network. Traffic is carried out in accordance with the assumptions of the \textbf{NaSch} model, with a low probability of random events and a maximum speed of 1. According to the recommendations of the model creators, updating position on the road should be asynchronous, so with each iteration of the program we draw a different order of road updates. Roads $1$--$4$ are roads with a horizontal direction and roads $5$-$8$ with a vertical direction. Drivers on even roads follow the natural turn, and drivers on odd   go in the opposite direction. By performing a single move for each street, the obtained results are placed in the appropriate positions on the network. The next stage is the analysis of behavior at intersections. The first 50 steps are skipped to allow the entire system to be filled with cars. The area before and at the intersection is taken into account. As per the authors’ recommendations, the intersection results are also updated in a random order. In the event of a meeting of two drivers at the intersection, the individual waiting time for each of them is set. It is calculated according to the values from Table \ref{tab1}, minus $1$, because the time needed to cross the intersection is not taken into account. Waiting time is then used in the previously mentioned algorithms for updating road positions. A vehicle ordered to wait cannot increase its speed until the designated number of steps has elapsed. Simulations were carried out for various driver relations on the road.

\subsection{Outcome of simulations. } 
\begin{figure}[tbh!]
		\centering
		\includegraphics[width=8.5cm]{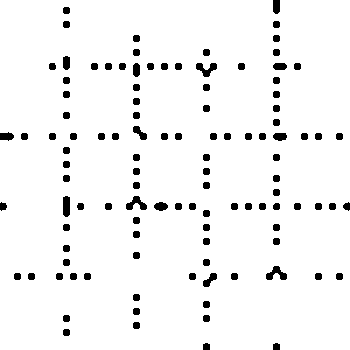}
		\caption{\label{siec1}The analyzed intersection network in one random program step.}
	\end{figure}

\begin{figure}[th!]
{\noindent
\small
\begin{subfigure}[b]{0.48\textwidth}
\centering {\includegraphics[width=0.8\textwidth]{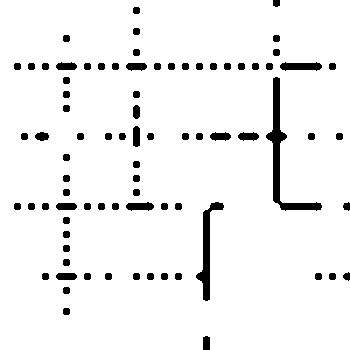}}
\caption{\label{siec2}The probability of occurrence of \textbf{DE} driver is 0.25.}
\end{subfigure}
\hfill\begin{subfigure}[b]{0.48\textwidth}
\centering {\includegraphics[width=0.8\textwidth]{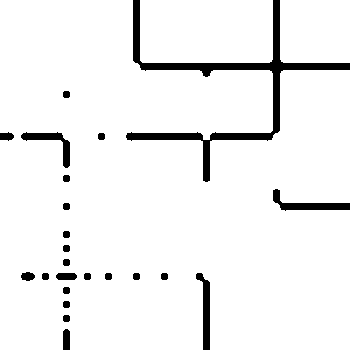}}
\caption{\label{siec3}The probability of occurrence of \textbf{DE} driver is 0.75.}
\end{subfigure}
}
\vspace{-1mm}
\caption{\label{siec2i3.}Analysis intersection networks in one random program iteration for more \textbf{DE} drivers.
}
\end{figure}

The parameters that can be modified are the maximum number of cars, the probability of introducing a new vehicle and the probability of occurrence of individual types of drivers. At the beginning, three examples are presented for different relations of drivers, for each the probability of a new car is $0.3$, and the maximum number of vehicles is $250$. A diagram will be presented from one selected moment (time step) of each simulation, the goal is to illustrate to the reader what the created network of intersections looks like and how the attitude of non-compliant drivers affects its traffic. The layout of the streets and the directions of vehicles follow the diagram in Fig. \ref{KWnetwork}. First, the network of intersections with smooth traffic was presented, for this purpose the probability of conflict driver will be $0.01$. The low value of this parameter means that such drivers hardly occur, hence, collision situations on the roads are rare. The described example is shown in Figure \ref{siec1}.
Another example (Fig. \ref{siec2}) is the situation of an increased number of non-compliant drivers. They occur with the  probability of $0.25$.  Traffic jams are noticeable, as shown by a number of cars waiting before the intersection. The reason for such events is the meeting of two conflict drivers, which leads to collisions and blocks the intersection.	
Severe traffic jams are presented in Fig. \ref{siec3}. The probability of a non-compliant driver is being raised again to $0.75$. At most intersections, these drivers meet, which causes conflicts and prevents cars from continuing to drive.	

\begin{figure}[tbh!]
		\centering
		\includegraphics[width=10cm]{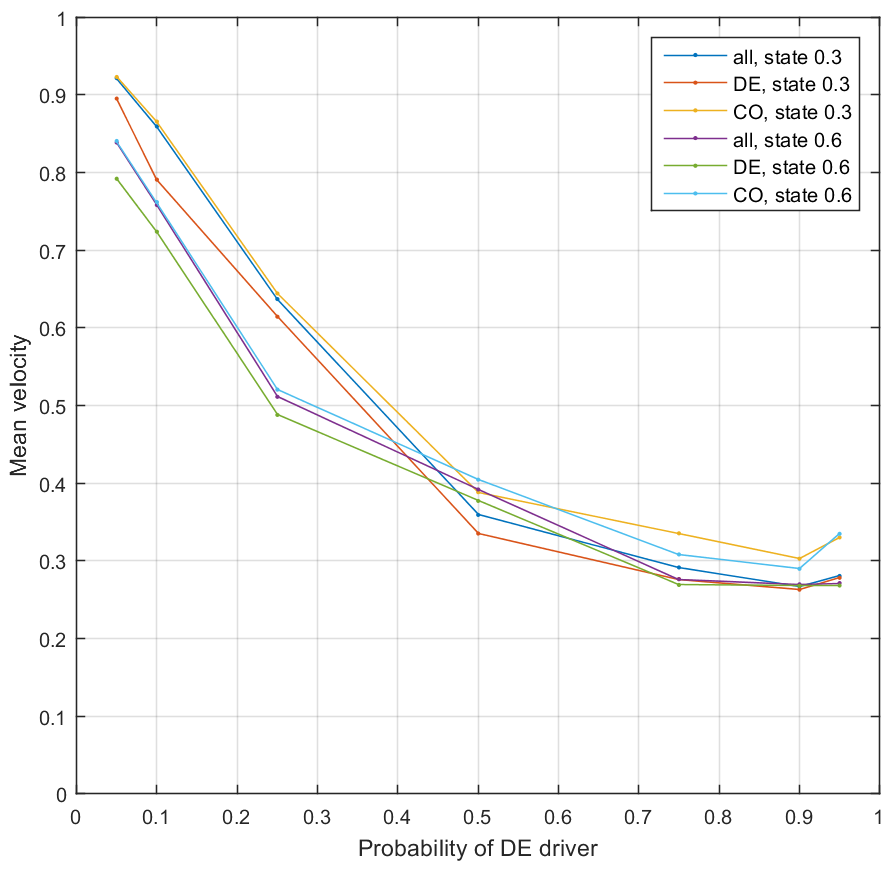}
		\caption{\label{razem}The average speeds of individual drivers type.} 
\end{figure}

In order to check how the occurrence of individual types of drivers affects the efficiency of the traffic network, average car speeds were compared depending on the ratio of drivers on the road. The probability of a non-compliant driver was increased, with the fixed probability of a new car of $0.3$ and a maximum number of drivers of $350$. At each step, the average speed of all vehicles as well as the average speed of each type of a driver was tested. $10,000$ replicates were carried out for each case, finally calculating the average of the values obtained. The process was repeated while increasing the probability of a new vehicle to $0.6$. Results are presented on picture \ref{razem}, here parameter \textit{state} indicates the probability of a new driver.
The first important conclusion is the fact that the average speeds of \textbf{CO} drivers -- that is, those who comply with traffic rules in each case are higher than for the other groups. It can be concluded that the movement of cooperating drivers is faster and smoother than the one of \textbf{DE} drivers who are more exposed to greater penalties when interacting at intersections. In addition, for high load of drivers on the road, with low probability of DE drivers, the average speeds are lower, but the presence of DE drivers improves traffic.	

\subsection{\label{KWKSzExtModInt}Extension of model interpretation.} The important issue for practice is \textbf{the estimation of driver type ratio}.
There are three types of meetings in the model under consideration: (\textbf{CO} vs \textbf{CO}; \textbf{CO} vs \textbf{DE}; \textbf{DE} vs. \textbf{DE}). Each of these types of meetings consequently gives one of the three effects of stream modification over the main road. Let $\eta_{ji}$ be a random variable equal to $1$ when in the $i$-th meeting of drivers, the driver on road $j$ will be \textbf{CO} and $0$ if he is \textbf{DE}, and $\xi_i = \eta_{ai} + \eta_{bi}$. If the random variables $\eta_{ji}$, $j \in \{a, b\}$ are independent, identically distributed with ${\bf P}(\eta_{ji} = 1) = p = 1-{\bf P} (\eta_{ji})$, then parameter $p$ can be estimated by one of the methods described in \cite{BicDok2015:MS} or  \cite{Ber1980:SDT}. 
Based on $n$  meetings of drivers with $\Sigma_{\xi} = \sum_1^n \xi_i$, the minimax estimate $\hat{p}_{minmax} = \frac {a + \Sigma_{\xi}}{b + 2n}$ (cf. \cite{MWal1971:SDF}, \cite{Ste1957:Estimation}).

Two cars are approaching an intersection ( or otherwise interacting due to e.g. a change of lane). Their drivers can follow the rules of the Highway Code or exceed the established rules (let's not figure out how). Despite the general rules of the road, there are situations in which the driver can choose his behavior to a certain extent. Some of them are fully compliant with traffic rules, and some are risky in the sense that they interfere with other users' roads to break their wrights. Sometimes the possible choices are limited by the behavior of other drivers. We can assume that this is the “determined property of the driver”, but also his conscious behavior - and thus the strategy. The first interpretation leads to the recognition that all road users are divided into “road users complying with traffic rules” or “those who violate these rules” This, in turn, leads to four types of meetings. It is seen that the proper \textbf{modeling drivers strategies} is crucial in the topic under consideration.

With the second interpretation, we can speak of a decision problem. Drivers do not have a permanently assigned feature, but only consciously make one of the two decisions. The mathematical model of this situation is a two-person game with a non-zero sum with a finite action space for both players (v. \cite{HarSel1988:GTgames} and \cite{Tij2003:Introduction} section 7, and 
\cite{Own2013:Game,Maz2014:book}). The players payoffs in this game are measured by the impact on their movement, mainly on speed. 

\section{\label{model2}Model II as an extension of the first  model.}
\unskip
\subsection{Description of the model.} 
The second model is an extension of the model proposed in section \ref{DetAnalysisModel1} with a constant probability of occurrence of a given type of a driver. The psychological model presenting the imitation strategy is taken into account (\cite[Chapter 15.3.1.]{Pla2012:Introduction}). The goal is to illustrate how the interaction between drivers affects the attitude of drivers on the road and, as a result, on the quality of traffic. The authors, inspired by earlier works from different areas, decided to apply a psychological model, where drivers follow the "do as others" principle. As before, two types of drivers were introduced, complying with and not complying with traffic rules. Drivers adopt strategies, not because of faith or a sense of duty to comply, but because of they imitate the behavior of others. In addition, a group of drivers who are not susceptible to the influence of other participants and faithfully following traffic rules is included. This group we will call \textbf{core}. The driver type is updated every $ \tau $ simulation steps. After this time, the probability of imitating each strategy is calculated, informing about the chance to change the current type of a driver to another with the probability of imitating it.

Such behavior is common in everyday life. In addition, many psychological works present such a model of society learning. Driver intelligence is not included in the participants' description road traffic. Many road users try to imitate others. This assumption is reflected in the intellectually underdeveloped environments. In addition, the core driver are included. It is unreasonable that drivers change their type of behavior too often, which is why determining the frequency of meeting individual types of drivers, based on which the driver changes his/her type of behavior, will be updated every $ \tau $ of simulation steps. On this basis, the imitation probability for each strategy will be calculated. It informs about the probability of changing the current type of a driver to another with the probability of imitation specified for him/her.  The type of driver observed is the result of his attribute and adopted strategy. This is sufficient for the purposes of this research to combine these into one parameter, although we expect interesting conclusions from the use of models based on hidden Markov chains (cf. \cite{LiHeZhou2016:HMM}, \cite{DenWuLyu2017:HMM}). 

More precisely, a change of a driver type from \textbf{CO} to \textbf{DE} (the road-complying driver to non-complying driver) occurs with the probability of $ P_{D} $, and in the reverse with the probability of $ P_{C} $. These probabilities are described by the formulas:
	\begin{equation} 
	P_{D} = \frac{f_{D}}{f_{C}+f_{D}}=1-P_{C},
	\label{pstwo:co}
	\end{equation} 
In the simulation $ f_s $ is the number of interactions of a given driver with the behavior of the opponent's type $s$ in a given measurement cycle $ \tau $. Traffic participants assess what driver their competitor was. The left driver, being a cooperator, cannot assess who his/her  opponent was. In such situations, we add $0.5$ to both $ f_{C} $ and $ f_{D} $. In other situations, we add $1$ to the appropriate counter. $\tau $ is set to 500.
	
\subsection{Simulation and results. } A simulation was carried out to check how the attitude of drivers changes in the imitation process depending on the size of the initial group of drivers who do not respect traffic rules. The effect will be observed by the ratio 
\begin{equation}\label{KWKSzratioDE}
\text{\cjRL{q}}=\frac{\#\textbf{DE}}{\#\textbf{CO}+\#\textbf{DE}}.
\end{equation}
\begin{figure}[h!]
{\noindent
\small
\begin{subfigure}[b]{0.95\textwidth}
\centering {\includegraphics[width=0.95\textwidth]{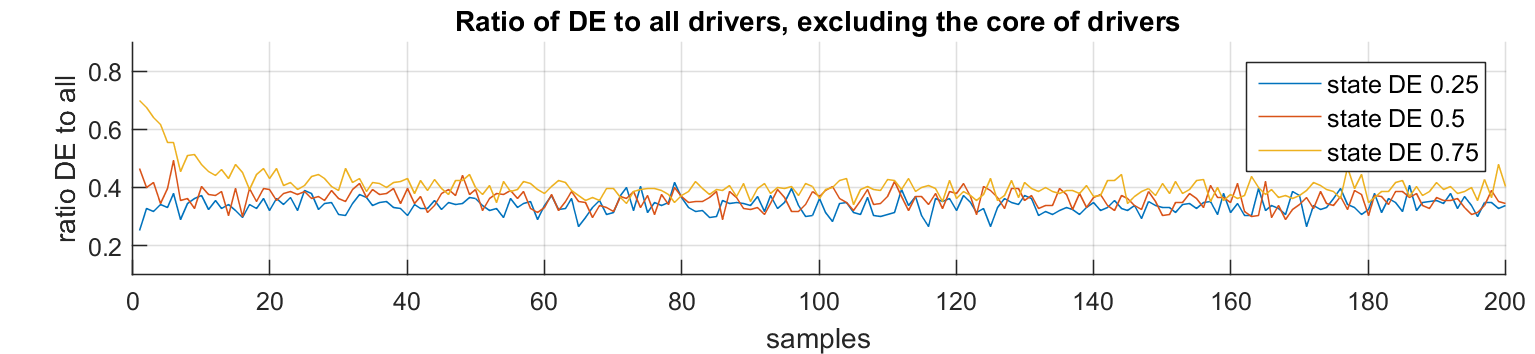}}
\caption{\label{psycho1}Excluding the core of drivers who comply with traffic rules.}
\end{subfigure}
\bigskip

\begin{subfigure}[b]{0.95\textwidth}
\centering {\includegraphics[width=0.95\textwidth]{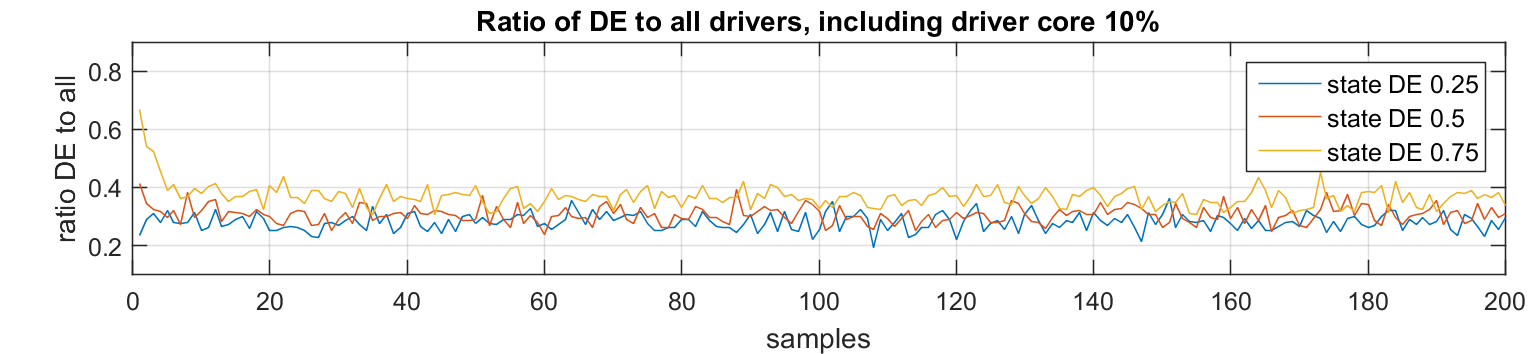}}
\caption{\label{psycho2} A driver core of 10\% of all drivers was included.}
\end{subfigure}
\bigskip

\begin{subfigure}[b]{0.95\textwidth}
\centering {\includegraphics[width=0.95\textwidth]{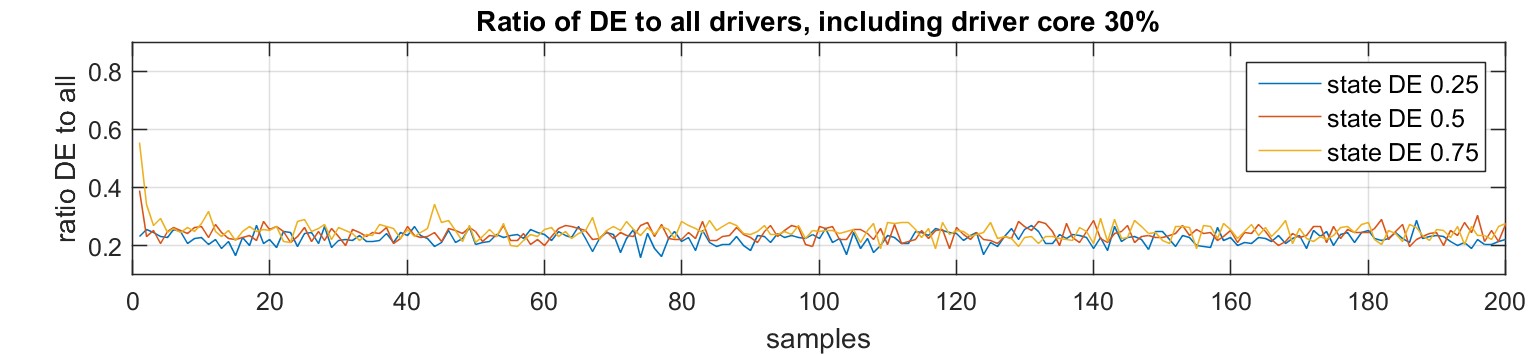}}
\caption{\label{psycho3} A driver core of 30\% of all drivers was included.}
\end{subfigure}
}
\vspace{-1mm}
\caption{\label{psycho123}The ratio \cjRL{q}. Different line colors correspond to the initial probabilities of \textbf{DE}. }
\end{figure}
Each simulation was carried out $ 200 \tau $ times, after each strategy update checking the ratio \text{\cjRL{q}}. Three models included. In the first, we assume that there is no permanent group (\textbf{core}) of drivers who comply with traffic rules and resist attempts to force them to break these rules. In the second and third models, we assume that there is a core of resistant and law-abiding drivers of 10\% and 30\% in total, respectively.

Figure \ref{psycho1} presents the results in the absence of drivers' core. The chart shows the ratio of \textbf{DE} drivers to all participants in situations where the initial probability of DE driver occurrence was $0.25$, $0.5$ and $0.75$, respectively, as shown by different line colors. It is noticeable that the values converge and remain at a similar level. It can be stated that the system is stabilizing in terms of the distribution of individual types of drivers. A similar situation is presented in Figure \ref{psycho2} showing simulation results with a driver core of 10\%. The system stabilizes at a lower level than before. In addition, the case with a high probability of occurrence of DE driver stands out more from the others. The situation after increasing the driver core to 30\% is presented in Figure \ref{psycho3}. It is important that this time the driver attitude stabilizes faster than in other cases, and the level of stabilization is even less than in the case of the core of 10\%. In addition, the charts for the different initial probability of \textbf{DE} are more similar. 
This presents an important fact resulting from the above analysis, the core of drivers significantly affects the level of stabilization of number of DE and CO drivers.  
The greater the core, the lower the stabilization level for DE drivers, which, when combined with the results of the previous model, 
gives a better flow and efficiency of the movement system. So, in order to ensure better quality of traffic, the emphasis should be on generating a larger core of drivers so that as many of them as possible are resistant to the negative influence of other road users.	
\begin{figure}[tbh!]
		\centering
		\includegraphics[width=14.5cm]{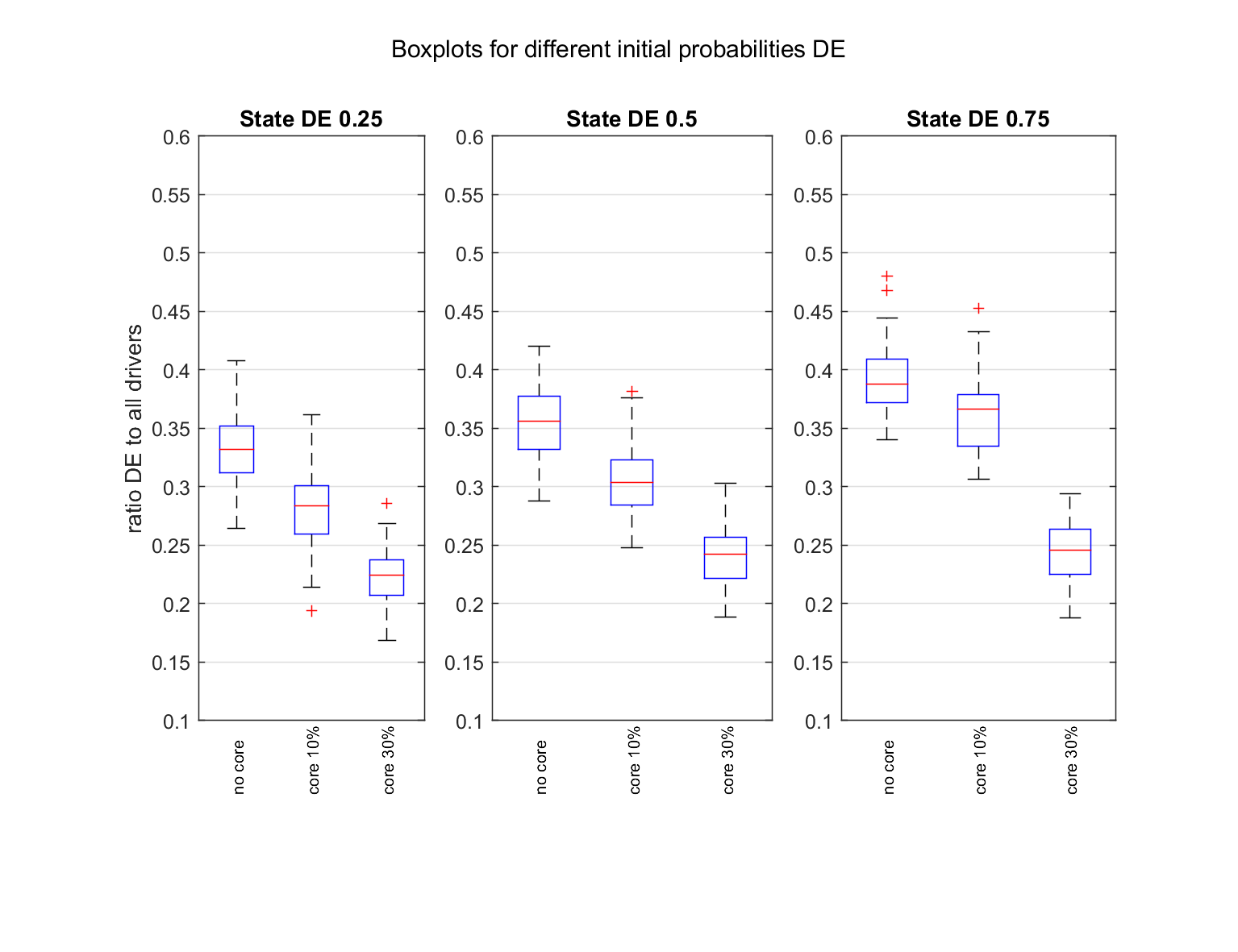}\\[-8ex]
		\caption{\label{psycho5}Box plots for the various initial probabilities of \textbf{DE} drivers vs. different \textbf{core} of \textbf{CO} drivers. 
		}
\end{figure}

The results obtained in the above simulations were summarized in box-and-whisker plots, previously removing $100$ initial values to ensure that the systems are in the stabilization phase. Values of different considerations of the driver core for each of the initial probabilities of DE drivers were compared. Figure \ref{psycho5} presents grouped charts for subsequent probabilities, and each of them has three box charts, presented in the following order:
\begin{enumerate}
		\item The case with a lack of the drivers' core complying with the regulations.
        \item The case where the drivers' core was 10\%.
        \item The case when the drivers' core was 30\%.
	\end{enumerate}
Presenting the results in this way confirms the earlier thesis that with the increase of the core of drivers the level of system stability is lower. In addition, with a larger initial ratio of DE drivers to the total, the size of the drivers' core has a greater impact on the final result, as seen in box-and-whisker 3 of Fig.~\ref{psycho5}. We note that the inter-quartile range for the results obtained is similar with maximum equal to about $0.4-0.5$. This confirms the fact that each of the analyzed systems is stabilized.
\begin{figure}[H]
\centering
	\includegraphics[width=14.5cm]{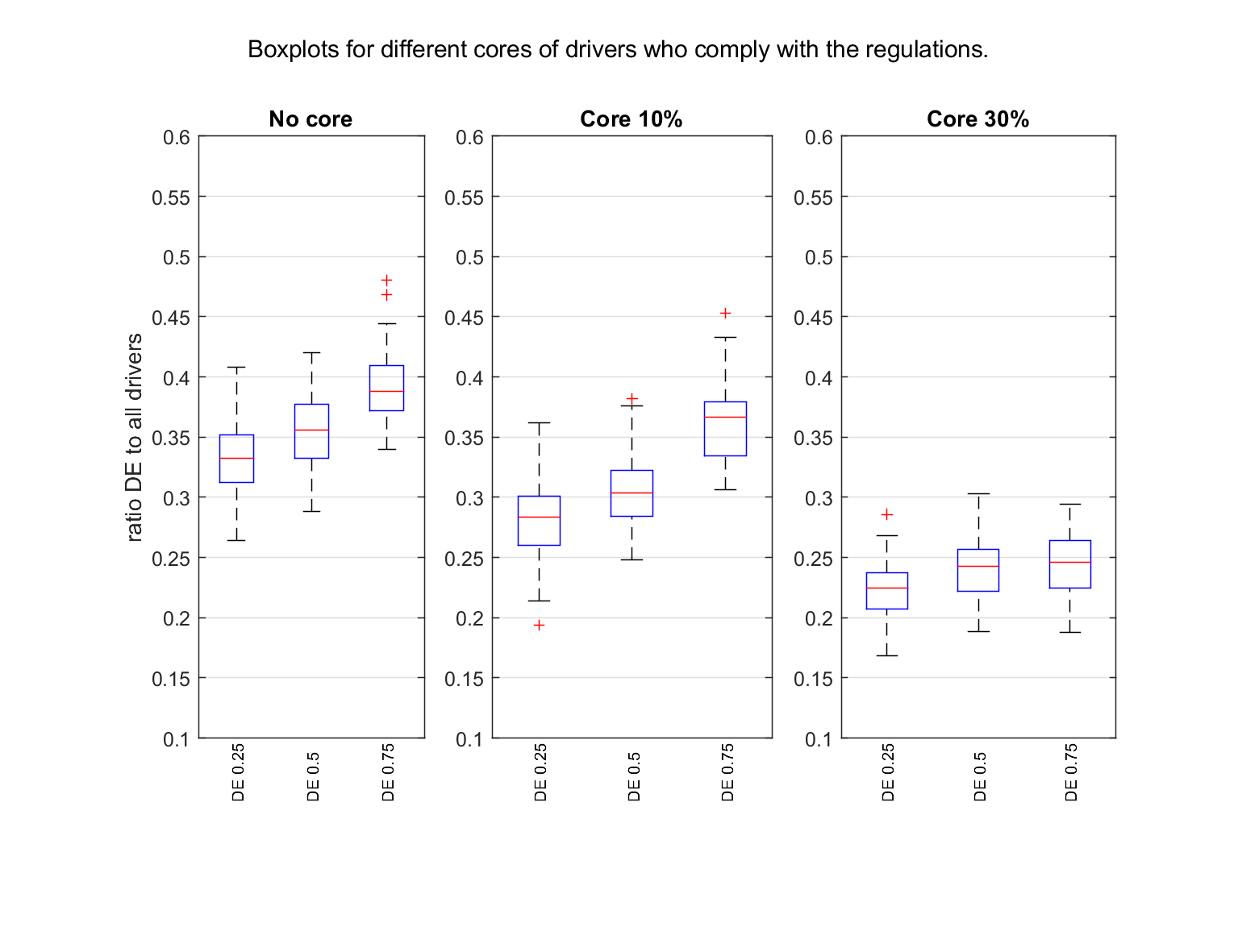}\\[-8ex]
		\caption{Box plots for various initial probabilities of \textbf{DE} drivers vs. the size of the \textbf{CO} drivers' \textbf{core}. 
		}\label{psycho6}
\end{figure}

The same box-and-whisker plots are compiled in a different way in Figure~\ref{psycho6}. In the previously presented three charts, this time drivers are divided by the size of the cores of prudent drivers. Each of them contains three box-plots, juxtaposed due to the initial probability of DE driver. This is a reflection of the graphs \ref{psycho1}--\ref{psycho3}, so the boxes successively indicate the initial probability of $0.25$, $0.5$ and $0.75$,  respectively. This confirms previous conclusions that the model with the driver's core of 10\% is the least stabilized, and the difference between the level of location of the chart for DE drivers with the initial ratio of $0.75$ is the largest. Hence, the model with a driver core of 30\% is the best stabilized. The final ratio of DE drivers to the whole is the lowest.

\section{\label{model3}Model III with impatient drivers.}
\unskip
\subsection{Problem formulation.} 
The last model considered was proposed in \cite{FanJiaTiaYun2014:TraficGT}. As in other models, the authors introduce two types of drivers (complying and not complying with the regulations). An important difference is how to generate individual types. It was assumed that at first everyone obeyed the rules of traffic, but after a certain time waiting before the intersection may cease to comply with the rules. This assumption is to reflect the actual behavior of drivers. Movement of vehicles, as in previous models, is simulated by cellular automata. The player's payouts, in this case are not explicitly stated, only the strategies that drivers use in each situation are known.

\subsection{\label{generating}Generating driver behavior.} As previously mentioned, the type of a driver depends on his/her  waiting time before the intersection. Drivers, waiting before the intersection, initially comply with the rules, but if the waiting time exceeds a certain individual value, the driver's behavior may change. This critical value is not constant and can be presented as the length of time the object was waiting to enter the intersection. It was assumed to be compatible with the Weibull distribution with the following cumulative distribution (cf. \cite{FanJiaTiaYun2014:TraficGT})
\begin{equation} 
	F(x) = \left\{ \begin{array}{ll}
	1-\exp \left\{ -\left( \frac{x}{a} \right)^b \right\}, &\text{for $x>0$,}\\
	0, & \textrm{otherwise,}\\
	\end{array} \right.
	\label{dystrybuanta}
\end{equation}
where $a$ is the scale parameter, $b$ is the shape parameter, $a,b$ are positive. The hazard rate function (on $\Re^+$) is
\begin{equation}  
	h(x) = \frac{f(x)}{1-F(x)}=\frac{b}{a} \left( \frac{x}{a} \right)^{b-1}. 
	\label{intensywnosc}
\end{equation} 
According to the proposed model, we set the scale parameter $a = 30$ and the shape parameter $b = 2.92$. These values represent the likelihood of changing driver behavior. When the driver begins to wait before the intersection, his/her behavior will change with the probability depending on the value of the function $ h(x) $ for a given waiting time. After passing an intersection, the driver's behavior returns to its initial state.	

\subsection{Driver strategies for prioritization.} Just as in the previous models, four types of interaction between drivers are possible: two drivers complying with the rules, two drivers not complying with the rules, two different drivers, the driver complying with the rules on a subordinate road and the driver not complying with the rules on a subordinate road. The following scenarios were highlighted:
\begin{enumerate}
\item The driver who should step down complies with the rules. By following them, it gives way regardless of the type of a driver at the cross-roads.
\item The driver who should give way is a driver who does not comply with the rules. He will try to impose priority, thus forcing the compliant driver to give way to him/her to ensure his/her safety.
\item Both drivers are non-compliant drivers, so both can try to cross the intersection at the same time. Because of their abnormal behavior, both of them should stop and then the driver on the road can pass first.
\end{enumerate}
	
\subsection{Simulation.} The assumptions of the above model have been implemented in the intersection network proposed in figure \ref{KWnetwork}. As in previous simulations, the maximum speed of cars was set at $1$ to explain the costs incurred by drivers when passing the intersection. In order to increase the likelihood of changing the type of behavior of a driver waiting before an intersection, a driver may enter the intersection when another leaves it. Table \ref{tab2} presents the time losses incurred by drivers while waiting before entering the intersection.
\begin{table}[h!]
		\begin{center} 
			\caption{Loss of time incurred by drivers during a meeting at the crossroads.}\label{tab2}
			\begin{tabular}{|c|c|c|} 
				\hline 
				{\backslashbox{Left-h. s.\strut}{\strut Right-h. s.}} & CO & DE\\
				\hline
				CO & 2,0 & 2,0\\
				\hline
				DE & 0,2 & 3,1\\
				\hline
			\end{tabular}
		\end{center}	
\end{table}
The left-hand side driver, when he is CO-type, always gives way. He waits two time steps, which is as much as his opponent needs to enter and leave the intersection. Then, if he is not moving, he can start the maneuver, and if another opponent arrives, the situation repeats. Waiting times before the intersection of each driver are counted. Additionally, drivers whose stopping does not result directly from waiting before the intersection, but is caused by its earlier blocking are also included. We consider cars standing in a traffic jam to those whose average of the previous five speeds is less than or equal to 0.2. Based on the received waiting times, the probability of change for each driver is determined.	

Let us analyze of the effect. Simulations were carried out for different probabilities of generating a new car, which is closely related to the density of cars on the road. A simulation was carried out for each case and was repeated 75,000 times. The following statistics were determined for each case:
\begin{enumerate}
\item Average system speed.
\item The number of a driver type changes that have occurred at each time step.
\item The number of conflicts between drivers that occurred at each time step.
\item The attitude of DE drivers, i.e. those who changed their type to non-compliant.
\item Average waiting times for drivers before intersections.	\end{enumerate}
The statistics listed are summarized in Table \ref{tab3}. It lists the number of conflicts and the number of all driver type changes that occurred in 75,000 repetitions for each generated case. The frequency of changes and the frequency of conflicts are considered in the analyzed samples. The average ratio of DE-type drivers to other drivers, showing the population of drivers who stopped complying with traffic rules, and average waiting time of drivers before the intersection.	
{	\small
	\begin{table}[th!]
		\begin{center} 
			\caption{\label{tab3}Summary of Model III Results}
			\begin{tabular}{|b{3.4cm}|c|c|c|c|c|c|c|c|c|}\hline 
				 & 0.1 & 0.2 & 0.3 & 0.4 & 0.5 & 0.6 & 0.7 & 0.8 & 0.9\\
				\hline
	{\small Total number of driver type changes} & 6 & 66 & 380 & 7,220 & 16,351 & 19,873 & 21,056 & 22,457 & 23,128\\
				\hline
	{\small Frequency of driver type changes} & 0.0001 & 0.0009 & 0.0051 & 0.0963 & 0.2180 & 0.2650 & 0.2807 & 0.2994 & 0.3084\\  
				\hline 
	{\small The total number of conflicts} & 2 & 0 & 0 & 25 & 218 & 248 & 376 & 337 & 388\\
				\hline\vspace{0.1cm}
	{\small Frequency of conflicts} & 0.00002 & 0 & 0 & 0.0003 & 0.0029 & 0.0033 & 0.0050 & 0.0045 & 0.0052\\
				\hline 
	{\small Average ratio of DE drivers} & 0.0073 & 0.0001 & 0.0004 & 0.0071 & 0.0197 & 0.0262 & 0.0184 & 0.0244 & 0.0231\\
				\hline
	{\small Average wait time} & 1.4527 & 1.3126 & 1.3759 & 2.4733 & 3.1016 & 3.3106 & 3.3708 & 3.4362 & 3.4878\\
				\hline
			\end{tabular}
		\end{center}	
	\end{table}
}

\begin{figure}[H]
	\centering
	\includegraphics[width=12cm]{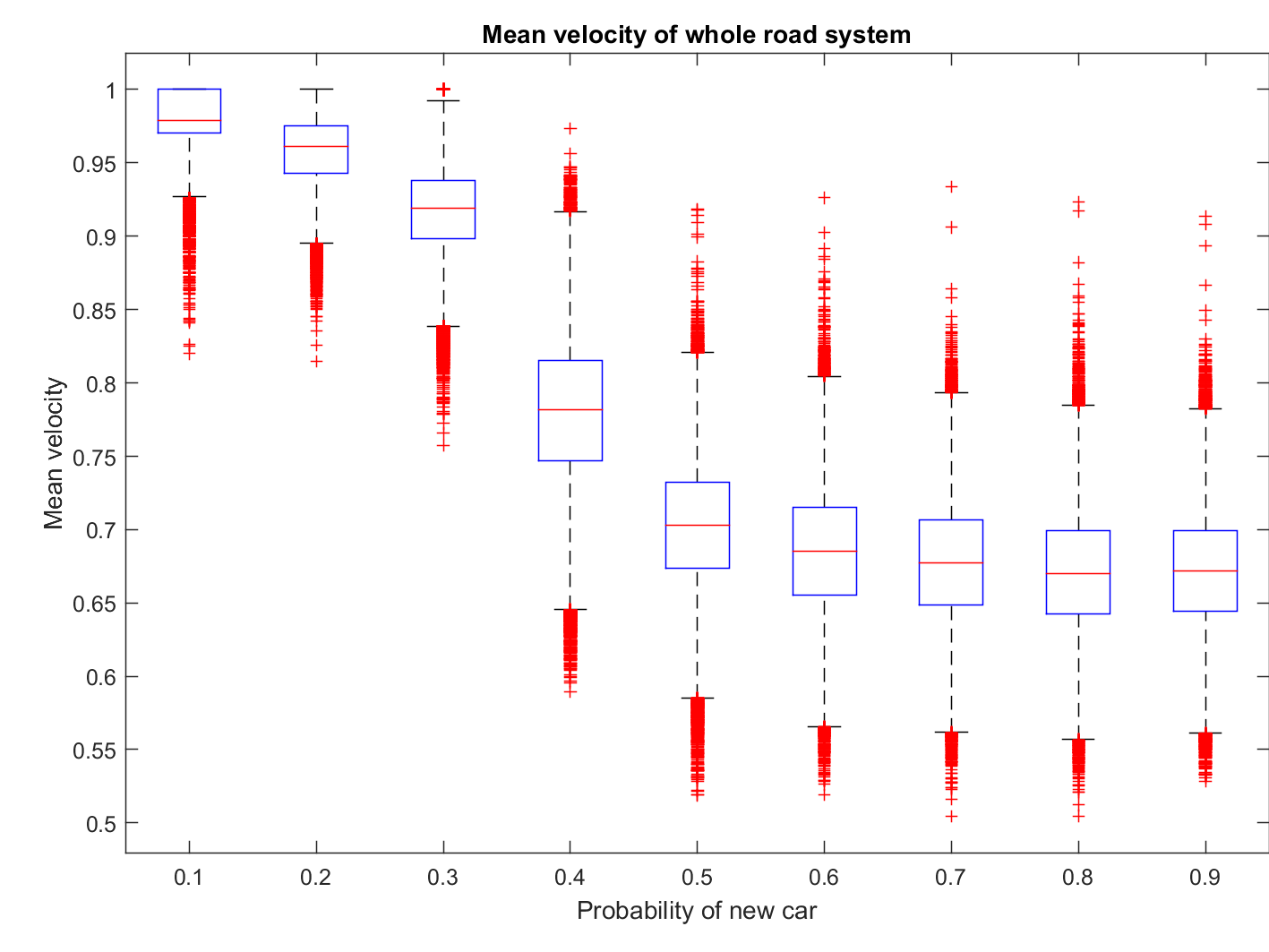}
	\caption{Box charts of average vehicle speeds depending on the likelihood of a new car occurring.}
	\label{new}
\end{figure}

The number of changes and the number of conflicts increase as the likelihood of a new driver increases. This is the result of high traffic density and traffic jams. Figure \ref{new} presents box charts of average speeds of the entire system from the entire motion simulation. A decrease in the efficiency of the traffic system is noticeable, the sharpest decrease occurs for the probability of a new driver equal to 0.4. This may represent a point where traffic density is too high and this causes traffic jams at intersections and waiting times before intersections are too long.

Based on the above results, it can be concluded that from the probability of drivers equal to 0.4, we are dealing with a high density of traffic with the efficiency of the system decreasing. So that the situation of too long waiting is highly likely, and drivers stop complying with traffic rules. The average waiting times in such situations are around $2.47$. Comparing this result with the values of the driver probability change function feature presented in Chapter \ref{generating}, 
it can be seen that these values are low. However, these are average values, at the same time one of the drivers could wait much longer, and another was just starting the waiting process. In addition, the driver, already standing in a traffic jam quite far from the intersection, can change his/her type, and return to the previous state only after leaving the intersection.

\section{\label{ModComp}Model Comparison.} 
The above work presents three ways of modeling traffic at intersections,focusing and the difference in modeling driver behavior. Each of them draws attention to a different problem regarding the functioning of traffic at intersections and each of them has been formulated in such a way that the obtained effects illustrate these problems and their consequences. The first presented model, which was dependent on the constant probability of receiving individual types of drivers, shows how the presence of non-compliant drivers negatively affects the quality of road traffic. It should be noted that the costs received by drivers in conflict situations are quite high, the purpose of this procedure is to draw attention to the negative effects of non-compliance with traffic rules. In fact, no collisions occur frequently. There are not so many non-compliant drivers on the roads and most road users are aware of the consequences of reckless driving.

However, the presented model makes it possible to illustrate the scale of the problem and possible effects if there were more reckless drivers. Another model, which is an extension of the previous one, draws attention to the psychological aspect of movement participants. It is assumed that drivers keep the imitation strategy, learn from each other regardless of the costs incurred. The effects of the received simulations show that regardless of the initial ratio of individual types of drivers, due to mutual learning, the distribution of individual categories of behavior converge to similar values. In addition, it is visible that, considering a certain group of drivers resistant to imitation strategy and always deciding to follow the law, we are able to influence the final distribution of types of behavior. This illustrates the potential of the correct education of future drivers. The last model can be considered a model of almost perfect traffic, because the drivers are mostly those who comply with the rules. In real life, traffic participants try to comply with traffic rules, however, there are situations when the patience of drivers reaches certain limits. During a long stand at an intersection, most drivers decide to violate the rules and force priority to be able to continue driving. It is natural that with more dense traffic this happens more often. It has been detected that above a certain degree of congestion, the traffic situation reaches certain limits of good performance. To sum up, each of the models in a different way reflects the processes occurring in the functioning of road traffic and draws attention to its different areas and problems. Therefore, conducting computer simulations allows to predict their effects.	

\section{\label{conclude}Conclusions.} 
The methodology of road traffic modeling connected with the drivers' behavior description is provided. The presented three driver behavior type models were combined with the Nagel-Schreckenber's models of car movements. This method of describing car movement is discussed in Section 1.3. Its link with driver types is presented in Section 2. Implementation of the models allowed their investigation based on the simulations. The maximum speeds equal to one were adopted, which, similarly, was applied in simple models. Their effects were checked on the basis of car traffic in a system consisting of 8 streets arranged perpendicularly, forming a network. As it was expected, various observations, qualitative and quantitative, were obtained. The applied approach allows forecasting car traffic based on knowledge of statistics of certain drivers' characteristics. To do this, you need to calibrate the model you want to use. The choice of the model depends on the possibility of obtaining information on the considered features of the driver population. Next to each of the considered models is a brief discussion on the possibilities of estimating model parameters from observing real traffic in the area under study. Thanks to this, the presented models can be used for traffic forecasting, calibration of driver behavior models and planning traffic protection by introducing its control.

They presented a different technique of linkage the drivers' behavior during interaction with others. The conclusions received were presented in the description of each simulation and the illustrated models were also compared. The conclusions drawn from the analyzed models affect the perception of various aspects of road traffic. For example, a large number of drivers who do not comply with traffic rules significantly worsens its quality and the presence of drivers who are faithful to their own views of compliance with traffic rules is able to improve its quality. Also, in the case of the model, when the driver’s behavior changes depending on the waiting time before the intersection, a relationship has been noticed between the density of cars on the road and the desire for drivers not to comply with the rules (v. Han and Ko~\cite{HanKo2012:junction}).


\vspace{6pt} 



\authorcontributions{Both authors equally contributed to the  conceptualization, methodology, formal analysis, investigation and writing--original draft preparation. Kinga Włodarczyk is responsible for the simulation software development, validation, visualisation and Krzysztof J. Szajowski is responsible for the   project administration and funds acquisition. 
}

\funding{This research received no external funding
}

\acknowledgments{The research is included in the leading topics of investigation of Faculty of Pure and Applied Mathematics, Wrocław University of Science and Technology under the project 049U/0051/19. 
}

\conflictsofinterest{The authors declare no conflict of interest.
} 

\abbreviations{The following abbreviations are used in this manuscript:\\

\noindent 
\begin{tabular}{@{}ll}
CO   						& \textbf{Regulatory drivers} (Cooperator)\\
DE   						& \textbf{Non-compliant drivers}(Defector)\\
\textbf{NaSch} 	& \textbf{Nagel-Schreckenberg}\\
CP, LCP         & Complementarity Problem, Linear Complementarity Problem\\
$\text{NE}(A,B)$& the set of Nash equilibria\\
\cjRL{q}        & The rate of non-compliant(jamming, difficult) drivers   
\end{tabular}}

	

\reftitle{References}


\externalbibliography{yes}

\bibliography{\jobname}

\begin{thebibliography}{-------}
\providecommand{\natexlab}[1]{#1}

\bibitem[Albi \em{et~al.}(2019)Albi, Bellomo, Fermo, Ha, Kim, Pareschi, Poyato,
  and Soler]{Alb2019:Vehicular}
Albi, G.; Bellomo, N.; Fermo, L.; Ha, S.Y.; Kim, J.; Pareschi, L.; Poyato, D.;
  Soler, J.
\newblock Vehicular traffic, crowds, and swarms: From kinetic theory and
  multiscale methods to applications and research
  perspectives${}^\textbf{\color{red}*}$.
\newblock {\em Mathematical Models and Methods in Applied Sciences} {\bf 2019},
  {\em 29},~1901--2005.
\newblock
  doi:{\changeurlcolor{black}\href{https://doi.org/10.1142/S0218202519500374}{\detokenize{10.1142/S0218202519500374}}}.

\bibitem[Useche \em{et~al.}(2017)Useche, Cendales, and
  Gómez]{UseCenGom2017:Work}
Useche, S.; Cendales, B.; Gómez, V.
\newblock Work stress, fatigue and risk behaviors at the wheel: Data to assess
  the association between psychosocial work factors and risky driving on Bus
  Rapid Transit drivers.
\newblock {\em Data in Brief} {\bf 2017}, {\em 15},~335 -- 339.
\newblock
  doi:{\changeurlcolor{black}\href{https://doi.org/10.1016/j.dib.2017.09.032}{\detokenize{10.1016/j.dib.2017.09.032}}}.

\bibitem[Chakravarty \em{et~al.}(2014)Chakravarty, Chowdhury, Ghose, Bhaumik,
  and Balamuralidhar]{ChaChow2014:Stat}
Chakravarty, T.; Chowdhury, A.; Ghose, A.; Bhaumik, C.; Balamuralidhar, P.
\newblock Statistical analysis of road-vehicle-driver interaction as an enabler
  to designing behavioural models.
\newblock {\em Journal of Physics: Conference Series} {\bf 2014}, {\em
  490},~012104.
\newblock
  doi:{\changeurlcolor{black}\href{https://doi.org/{10.1088/1742-6596/490/1/012104}}{\detokenize{{10.1088/1742-6596/490/1/012104}}}}.

\bibitem[Greenshields(1935)]{Gre1935:Traffic}
Greenshields, B.
\newblock A study of traffic capacity.
\newblock {\em Proceedings of the Highway Research Board} {\bf 1935}, {\em
  14},~448--477.

\bibitem[Zhang(1998)]{Zha1998:nonequilibrium}
Zhang, H.
\newblock A theory of nonequilibrium traffic flow.
\newblock {\em Transportation Research Part B: Methodological} {\bf 1998}, {\em
  32},~485--498.
\newblock
  doi:{\changeurlcolor{black}\href{https://doi.org/10.1016/S0191-2615(98)00014-9}{\detokenize{10.1016/S0191-2615(98)00014-9}}}.

\bibitem[Helbing \em{et~al.}(2001)Helbing, Hennecke, Shvetsov, and
  Treiber]{gass2001:Master}
Helbing, D.; Hennecke, A.; Shvetsov, V.; Treiber, M.
\newblock MASTER: macroscopic traffic simulation based on a gas-kinetic,
  non-local traffic model.
\newblock {\em Transportation Research B} {\bf 2001}, {\em 32},~183--211.
\newblock
  doi:{\changeurlcolor{black}\href{https://doi.org/10.1016/S0191-2615(99)00047-8}{\detokenize{10.1016/S0191-2615(99)00047-8}}}.

\bibitem[Liu \em{et~al.}(2007)Liu, Xin, Adams, and Ban]{LiuXinAdaBan2007:Game}
Liu, H.X.; Xin, W.; Adams, Z.M.; Ban, J.X.
\newblock A Game Theoretical Approach for Modelling Merging and Yielding
  Behavior at Freeway On-Ramp Sections.
\newblock  Transportation and Traffic Theory 2007; Taylor, M., Ed.; Elsevier.
  Linacre House.: Jordan Hill, Oxford, United Kingdom OX2 8DP,  2007; pp.
  196--211.
\newblock Transportation and Traffic Theory 2007. London, U.K. 23-25.07.2007.

\bibitem[Kita \em{et~al.}(2002)Kita, Tanimoto, and Fukuyama]{KitKei2002:Game}
Kita, H.; Tanimoto, K.; Fukuyama, K.
\newblock {A game theoretic analysis of merging -- giveaway interaction: A
  joint estimation model}. In {\em Transportation and Traffic Theory in the
  21st Century}; Taylor, M., Ed.; Emerald Group Publishing Limited: Adelaide,
  Australia,  2002; pp. 503--518.
\newblock
  doi:{\changeurlcolor{black}\href{https://doi.org/10.1108/9780585474601-025}{\detokenize{10.1108/9780585474601-025}}}.

\bibitem[Yao \em{et~al.}(2018)Yao, Jia, Zhong, and
  Li]{YaoJiaZhoLi2018:BestResponse}
Yao, W.; Jia, N.; Zhong, S.; Li, L.
\newblock Best response game of traffic on road network of non-signalized
  intersections.
\newblock {\em Physica A: Statistical Mechanics and its Applications} {\bf
  2018}, {\em 490},~386 -- 401.
\newblock
  doi:{\changeurlcolor{black}\href{https://doi.org/https://doi.org/10.1016/j.physa.2017.08.032}{\detokenize{https://doi.org/10.1016/j.physa.2017.08.032}}}.

\bibitem[Ilachinski(2001)]{Ila2001:Automata}
Ilachinski, A.
\newblock {\em Cellular automata}; World Scientific Publishing Co., Inc., River
  Edge, NJ,  2001; pp. xxxii+808.
\newblock A discrete universe,
  doi:{\changeurlcolor{black}\href{https://doi.org/10.1142/4702}{\detokenize{10.1142/4702}}}.

\bibitem[Neumann(1966)]{Neu1966:Self}
Neumann, J.V.
\newblock {\em Theory of Self-Reproducing Automata}; University of Illinois
  Press: Champaign, IL, USA,  1966.
\newblock Edited and completyt by Burks, Arthur W.

\bibitem[Emb(2010)]{Embryo2010}
{John von Neumann's Cellular Automata},  2010.
\newblock \href{http://embryo.asu.edu/handle/10776/2009}{Embryo Project
  Encyclopedia (2010-06-14)}.

\bibitem[Małecki and Szmajdziński(2013)]{MalSzm2013:AnalizaRuchu}
Małecki, K.; Szmajdziński, K.
\newblock Simulator for microscopic traffic analysis.
\newblock {\em Logistyka} {\bf 2013}, {\em 3},~8.
\newblock Title of Polish org. “Symulator do mikroskopowej analizy ruchu
  drogowego”.

\bibitem[Berto and Tagliabue(2017)]{BerTag2017:CAmath}
Berto, F.; Tagliabue, J.
\newblock
  \href{https://plato.stanford.edu/archives/fall2017/entries/cellular-automata/}{Cellular
  Automata}. In {\em The Stanford Encyclopedia of Philosophy}, Fall 2017 ed.;
  Zalta, E.N., Ed.; Metaphysics Research Lab, Stanford University,  2017.

\bibitem[Żygierewicz(2019)]{Zyg2019:AuKo}
Żygierewicz, J.
\newblock
  \href{http://www.fuw.edu.pl/~jarekz/MODELOWANIE/Automaty_komorkowe.pdf}{Cellular
  automata},  2019.
\newblock Title of Polish org. “Automaty komórkowe”. The Author’s web
  page: Home web page,
  \href{http://www.fuw.edu.pl/~jarekz}{http://www.fuw.edu.pl/~jarekz}, access:
  08.01.2019.

\bibitem[Tanimoto \em{et~al.}(2014)Tanimoto, Kukida, and
  Hagishima]{TanKukHag2014:Social}
Tanimoto, J.; Kukida, S.; Hagishima, A.
\newblock Social dilemma structures hidden behind traffic flow with lane
  changes${}^\textbf{\color{red}*}$.
\newblock {\em Journal of Statistical Mechanics: Theory and Experiment} {\bf
  2014}, {\em 2014},~P07019.
\newblock
  doi:{\changeurlcolor{black}\href{https://doi.org/10.1088/1742-5468/2014/07/P07019}{\detokenize{10.1088/1742-5468/2014/07/P07019}}}.

\bibitem[Fisk(1984)]{Fis1984:Game}
Fisk, C.S.
\newblock Game theory and transportation systems modelling.
\newblock {\em Transportation Res. Part B} {\bf 1984}, {\em 18},~301--313.
\newblock
  doi:{\changeurlcolor{black}\href{https://doi.org/10.1016/0191-2615(84)90013-4}{\detokenize{10.1016/0191-2615(84)90013-4}}}.

\bibitem[Ferguson(1967)]{Fer1967:MS}
Ferguson, T.S.
\newblock {\em Mathematical statistics: {A} decision theoretic approach};
  Probability and Mathematical Statistics, Vol. 1, Academic Press, New
  York-London,  1967; pp. xi+396.

\bibitem[Owen(2013)]{Own2013:Game}
Owen, G.
\newblock {\em Game theory}, fourth ed.; Emerald Group Publishing Limited,
  Bingley,  2013; pp. viii+451.
\newblock Tłum. polskie (tłum. Andrzej Wieczorek): PWN, Warszawa 1975.
  \MR{3443071}.

\bibitem[Płatkowski(2012)]{Pla2012:Introduction}
Płatkowski, T.
\newblock {\em
  \href{http://mst.mimuw.edu.pl/wyklady/wtg/wyklad.pdf}{Introduction to Game
  Theory}}; Warsaw University,  2012; p. 104.
\newblock Title of Polish org. "Wstęp do Teorii Gier". The Author's web page:
  \href{http://mst.mimuw.edu.pl/}{http://mst.mimuw.edu.pl/. Access:
  08.01.2019r.}

\bibitem[Mazalov(2014)]{Maz2014:book}
Mazalov, V.V.
\newblock {\em Mathematical Game Theory and Applications}; John Wiley \& Sons,
  Ltd., Chichester,  2014; pp. xiv+414.

\bibitem[Hidas(2002)]{Hid2002:Merging}
Hidas, P.
\newblock Modelling Lane Changing and Merging in Microscopic Traffic
  Simulation.
\newblock {\em Transportation Res. Part C} {\bf 2002}, {\em 10},~351--371.
\newblock
  doi:{\changeurlcolor{black}\href{https://doi.org/10.1016/S0968-090X(02)00026-8}{\detokenize{10.1016/S0968-090X(02)00026-8}}}.

\bibitem[Paissan and Abramson(2013)]{PasGui2013:Imitation}
Paissan, G.; Abramson, G.
\newblock Imitation dynamics in a game of traffic.
\newblock {\em The European Physical Journal B} {\bf 2013}, {\em 86},~153.
\newblock
  doi:{\changeurlcolor{black}\href{https://doi.org/10.1140/epjb/e2013-30372-5}{\detokenize{10.1140/epjb/e2013-30372-5}}}.

\bibitem[Fan \em{et~al.}(2014)Fan, Jia, Tian, and
  Yun]{FanJiaTiaYun2014:TraficGT}
Fan, H.; Jia, B.; Tian, J.; Yun, L.
\newblock Characteristics of traffic flow at a non-signalized intersection in
  the framework of game theory.
\newblock {\em Physica A: Statistical Mechanics and its Applications} {\bf
  2014}, {\em 415},~172 -- 180.
\newblock
  doi:{\changeurlcolor{black}\href{https://doi.org/https://doi.org/10.1016/j.physa.2014.07.031}{\detokenize{https://doi.org/10.1016/j.physa.2014.07.031}}}.

\bibitem[Ding and Huang(2009)]{DinHua2009:TrafficFlow}
Ding, Z.; Huang, G.
\newblock Real-Time Traffic Flow Statistical Analysis Based on
  Network-Constrained Moving Object Trajectories.
\newblock  Database and Expert Systems Applications; Bhowmick, S.S.; K{\"u}ng,
  J.; Wagner, R., Eds.; Springer Berlin Heidelberg: Berlin, Heidelberg,  2009;
  pp. 173--183.
\newblock
  doi:{\changeurlcolor{black}\href{https://doi.org/10.1007/978-3-642-03573-9\_14}{\detokenize{10.1007/978-3-642-03573-9\_14}}}.

\bibitem[Bifulco \em{et~al.}(2014)Bifulco, Galante, Pariota, Spena, and
  Gais]{BifGalParSpeGai2014:TrafficData}
Bifulco, G.; Galante, F.; Pariota, L.; Spena, M.R.; Gais, P.D.
\newblock Data Collection for Traffic and Drivers’ Behaviour Studies: A
  Large-scale Survey.
\newblock {\em Procedia - Social and Behavioral Sciences} {\bf 2014}, {\em
  111},~721 -- 730.
\newblock Transportation: Can we do more with less resources? – 16th Meeting
  of the Euro Working Group on Transportation – Porto 2013,
  doi:{\changeurlcolor{black}\href{https://doi.org/https://doi.org/10.1016/j.sbspro.2014.01.106}{\detokenize{https://doi.org/10.1016/j.sbspro.2014.01.106}}}.

\bibitem[Bickel and Doksum(2015)]{BicDok2015:MS}
Bickel, P.J.; Doksum, K.A.
\newblock {\em Mathematical statistics---basic ideas and selected topics.
  {V}ol. 1}, second ed.; Vol.~1, {\em Texts in Statistical Science Series}, CRC
  Press, Boca Raton, FL,  2015; pp. xxiv+556.

\bibitem[Berger(1980)]{Ber1980:SDT}
Berger, J.O.
\newblock {\em Statistical decision theory: foundations, concepts, and
  methods}; Springer-Verlag, New York-Heidelberg,  1980; pp. xv+425.
\newblock Springer Series in Statistics.

\bibitem[{Federal Highway Administration}(2009)]{MUTCD2009}
{Federal Highway Administration}.
\newblock {\em Manual on Uniform Traffic Control Devices (MUTCD)}; Vol. Part 1,
  United States Department of Transportation,  2009.
\newblock Access November 28, 2011.

\bibitem[Nagel and Schreckenberg(1992)]{NagSch1992:freeway}
Nagel, K.; Schreckenberg, M.
\newblock A cellular automaton model for freeway traffic.
\newblock {\em Journal de Physique I France} {\bf 1992}, {\em 2},~2221--2229.
\newblock
  doi:{\changeurlcolor{black}\href{https://doi.org/10.1051/jp1:1992277}{\detokenize{10.1051/jp1:1992277}}}.

\bibitem[Cortés-Berrueco \em{et~al.}(2016)Cortés-Berrueco, Gershenson, and
  Stephens]{Cor-BerGerSte2016:Traffic}
Cortés-Berrueco, L.E.; Gershenson, C.; Stephens, C.R.
\newblock Traffic Games: Modeling Freeway Traffic with Game
  Theory${}^\textbf{\color{red}*}$.
\newblock {\em PLOS ONE} {\bf 2016}, {\em 11},~1--34.
\newblock
  doi:{\changeurlcolor{black}\href{https://doi.org/10.1371/journal.pone.0165381}{\detokenize{10.1371/journal.pone.0165381}}}.

\bibitem[Mesterton-Gibbons(1990)]{MesGub1990:Dilemma}
Mesterton-Gibbons, M.
\newblock A game-theoretic analysis of a motorist's
  dilemma${}^\textbf{\color{red}*}$.
\newblock {\em Math. Comput. Modelling} {\bf 1990}, {\em 13},~9--14.
\newblock
  doi:{\changeurlcolor{black}\href{https://doi.org/10.1016/0895-7177(90)90028-L}{\detokenize{10.1016/0895-7177(90)90028-L}}}.

\bibitem[Harsanyi and Selten(1988)]{HarSel1988:GTgames}
Harsanyi, J.C.; Selten, R.
\newblock {\em A general theory of equilibrium selection in games}; MIT Press,
  Cambridge, MA,  1988; pp. xviii+378.
\newblock With a foreword by Robert Aumann. \MR{956053}.

\bibitem[Wald(1971)]{MWal1971:SDF}
Wald, A.
\newblock {\em {Statistical Decision Functions}}; Chelsea Publishing Co.,
  Bronx, N.Y.,  1971; pp. ix+179.
\newblock Reprint of the 1950 edition. \MR{0394957}.

\bibitem[Steinhaus(1957)]{Ste1957:Estimation}
Steinhaus, H.
\newblock The problem of estimation.
\newblock {\em Ann. Math. Statist.} {\bf 1957}, {\em 28},~633--648.
\newblock
  doi:{\changeurlcolor{black}\href{https://doi.org/10.1214/aoms/1177706876}{\detokenize{10.1214/aoms/1177706876}}}.

\bibitem[Tijs(2003)]{Tij2003:Introduction}
Tijs, S.
\newblock {\em Introduction to game theory}; Vol.~23, {\em Texts and Readings
  in Mathematics}, Hindustan Book Agency, New Delhi,  2003; pp. viii+176.

\bibitem[Li \em{et~al.}(2016)Li, He, Zhou, Guan, and Dai]{LiHeZhou2016:HMM}
Li, J.; He, Q.; Zhou, H.; Guan, Y.; Dai, W.
\newblock Modeling Driver Behavior near Intersections in Hidden Markov Model.
\newblock {\em Int J Environ Res Public Health} {\bf 2016}, {\em 13},~15 pages.
\newblock \PMC{5201406};\PMID{28009838},
  doi:{\changeurlcolor{black}\href{https://doi.org/10.3390/ijerph13121265}{\detokenize{10.3390/ijerph13121265}}}.

\bibitem[Deng \em{et~al.}(2017)Deng, Wu, Lyu, and Huang]{DenWuLyu2017:HMM}
Deng, C.; Wu, C.; Lyu, N.; Huang, Z.
\newblock Driving style recognition method using braking characteristics based
  on hidden Markov model.
\newblock {\em PLOS ONE} {\bf 2017}, {\em 12},~1--15.
\newblock
  doi:{\changeurlcolor{black}\href{https://doi.org/10.1371/journal.pone.0182419}{\detokenize{10.1371/journal.pone.0182419}}}.

\bibitem[{Han, {Yo-Sub} and Ko, {Sang-Ki}}(2012)]{HanKo2012:junction}
{Han, {Yo-Sub} and Ko, {Sang-Ki}}.
\newblock Analysis of a cellular automaton model for car traffic with a
  junction.
\newblock {\em Theoret. Comput. Sci.} {\bf 2012}, {\em 450},~54--67.
\newblock
  doi:{\changeurlcolor{black}\href{https://doi.org/10.1016/j.tcs.2012.04.027}{\detokenize{10.1016/j.tcs.2012.04.027}}}.

\end{thebibliography}

\end{document}